\documentclass{svmult}
\usepackage{subeqn}
\usepackage[pdftex]{graphicx}
\usepackage{pstricks}
\usepackage{amsfonts,amssymb,slashbox}
\usepackage{mathptmx,helvet,courier,makeidx,multicol,footmisc}
\usepackage{natbib}
\bibpunct{(}{)}{;}{a}{}{,}
\usepackage[bookmarksnumbered=true, pdfauthor={Wen-Long Jin}]{hyperref}
\newcommand{\commentout}[1]{}

\newcommand{\ba}{\begin{array}}
        \newcommand{\ea}{\end{array}}
\newcommand{\bc}{\begin{center}}
        \newcommand{\ec}{\end{center}}
\newcommand{\bdm}{\begin{displaymath}}
        \newcommand{\edm}{\end{displaymath}}
\newcommand{\bds} {\begin{description}}
        \newcommand{\eds} {\end{description}}
\newcommand{\ben}{\begin{enumerate}}
        \newcommand{\een}{\end{enumerate}}
\newcommand{\beq}{\begin{equation}}
        \newcommand{\eeq}{\end{equation}}
\newcommand{\bfg} {\begin{figure}[h]}
        \newcommand{\efg} {\end{figure}}
\newcommand{\bi} {\begin {itemize}}
        \newcommand{\ei} {\end {itemize}}
\newcommand{\bqn}{\begin{eqnarray}}
        \newcommand{\eqn}{\end{eqnarray}}
\newcommand{\bqs}{\begin{eqnarray*}}
        \newcommand{\eqs}{\end{eqnarray*}}
\newcommand{\bsl} {\begin{slide}[8.8in,6.7in]}
        \newcommand{\esl} {\end{slide}}
\newcommand{\bss} {\begin{slide*}[9.3in,6.7in]}
        \newcommand{\ess} {\end{slide*}}
\newcommand{\btb} {\begin {table}}
        \newcommand{\etb} {\end {table}}
\newcommand{\bvb}{\begin{verbatim}}
        \newcommand{\evb}{\end{verbatim}}

\newcommand{\m}{\mbox}

\newcommand {\pd}[2] {{\frac {\partial {#1}} {\partial {#2}}}}
\newcommand {\pdd}[2] {{\frac {\partial^2 {#1}} {\partial {#2}^2}}}

\newcommand{\cas}[1]{{{\left \{ \ba #1 \ea \right. }}}

\newcommand{\reff}[1] {{{Figure \ref {#1}}}}
\newcommand{\refe}[1] {{{(\ref {#1})}}}
\newcommand{\reft}[1] {{{\textbf{Table} \ref {#1}}}}

\def\pmb#1{\setbox0=\hbox{$#1$}%
   \kern-.025em\copy0\kern-\wd0
   \kern.05em\copy0\kern-\wd0
   \kern-.025em\raise.0433em\box0 }

\def\eop{{\hfill $\blacksquare$}}
\def\r{{\rho}}

\def\dx     {{\Delta x}}
\def\dt     {{\Delta t}}

\begin {document}
\title*{Supply-demand diagrams and a new framework for analyzing the inhomogeneous Lighthill-Whitham-Richards model}
\author{Wen-Long Jin, University of California, Irvine, CA; Liang Chen, University of Science and Technology of China, Hefei, China; Elbridge Gerry Puckett, University of California, Davis, CA}
\titlerunning{Supply-demand diagrams and the inhomogeneous LWR model}
\authorrunning{Jin, Chen and Puckett}

\maketitle

\abstract{Traditionally, the Lighthill-Whitham-Richards (LWR) models for homogeneous and inhomogeneous roads have been analyzed in flux-density space with the fundamental diagram of the flux-density relation. In this paper, we present a new framework for analyzing the LWR model, especially the Riemann problem at a linear boundary in which the upstream and downstream links are homogeneous and initially carry uniform traffic. We first review the definitions of local supply and demand functions and then introduce the so-called supply-demand diagram, on which a traffic state can be represented by its supply and demand, rather than as density and flux as on a fundamental diagram.
It is well-known that the solutions to the Riemann problem at each link are self-similar with a stationary state, and that the wave on the link is determined by the stationary state and the initial state. In our new framework, there can also exist an interior state next to the linear boundary on each link, which takes infinitesimal space, and admissible conditions for the upstream and downstream stationary and interior states can be derived in supply-demand space. With an entropy condition consistent with a local supply-demand method in interior states, we show that the stationary states exist and are unique within the solution framework. We also develop a graphical scheme for solving the Riemann problem, and the results are shown to be consistent with those in the literature.
We further discuss asymptotic stationary states on an inhomogeneous ring road with arbitrary initial conditions and demonstrate the existence of interior states with a numerical example.
The framework developed in this study is simpler than existing ones and can be extended for analyzing the traffic dynamics in general road networks.}
\section{Introduction}
Essential to effective and efficient transportation control, management, and planning strategies is a better understanding of the evolution of traffic dynamics on a road network; i.e., the formation, propagation, and dissipation of traffic queues. The seminal work by \citep{lighthill1955lwr} and \citep{richards1956lwr} (LWR) attempts to study traffic dynamics with respect to aggregate values such as density $\r$, speed $v$, and flux $q$. Based on a continuous version of traffic conservation
 \bqn
 \pd{\r}{t}+\pd{q}x=0
 \eqn
 and an assumption about the fundamental diagram of the flux-density relation $q=Q(\r)$, the LWR model of a homogeneous road link can be written as
\bqn
\pd{\r}{t}+\pd {Q(\r)}x=0. \label{lwr}
\eqn
The corresponding speed-density relation is $v=V(\r)\equiv Q(\r)/\r$.
Here the maximum or jam density is denoted by $\r_j$; i.e., $\r \in[0,\r_j]$. Usually, $V(\r)$ is a non-increasing function of traffic density, $v_f=V(0)$ is the free flow speed, $V(\r_j)=0$, and $q=Q(\r)$ is unimodal with maximum flux or capacity $C=Q(\r_c)$ where $\r_c$ is the critical density.
Finally, traffic states with density higher than $\r_c$ are congested or over-critical, and those with density lower than $\r_c$ are free flowing or under-critical.

Compared with microscopic traffic flow models \citep[e.g.][]{gazis1961follow,nagel1992ca} the LWR model can be used to analyze traffic evolution at the aggregate level with shock and rarefaction waves. With its analytical power and simplicity, the LWR theory has been extended for studying traffic dynamics in more general transportation networks. For examples, \citet{daganzo1997special} proposed a traffic flow model for freeways with special lanes and high-occupancy vehicles with a two-regime fundamental diagram, and \citet{wong2002multiclass} proposed a multi-class model for heterogeneous drivers.

In this paper, we are interested in the LWR model for a road with bottlenecks, where traffic characteristics such as free flow speed, jam density, the number of lanes, and capacity may be different for different locations. In other words, the fundamental diagram $q=Q(x,\r)$ depends on location. Such a road link is called inhomogeneous and the corresponding inhomogeneous LWR model can be written as
\bqn
\pd{\r}{t}+\pd {Q(x,\r)}x=0. \label{inhlwr}
\eqn
In order to understand the fundamental properties of equation \refe{inhlwr}, we usually analyze its Riemann problem at $x=0$. Hereafter, we will refer to the upstream branch as link 1, the downstream branch as link 2, and $x=0$ as a linear boundary. In the Riemann problem, links 1 and 2 are both homogeneous and initially carry uniform traffic. That is,
\bqn
Q(x,\r)&=&\cas{{ll} Q_1(\r), & x<0, \\ Q_2(\r), &x>0,} \label{discontinuousflux}
\eqn
and
\bqn
\r(x,t=0)&=&\cas{{ll} \r_1, &x<0, \\ \r_2, &x>0.} \label{inicon}
\eqn

Since \citep{mochon1987changing}, there have been many analytical and numerical studies related to the inhomogeneous LWR model in the literature. Roughly speaking, there have been two types of methods for solving the Riemann problem of inhomogeneous LWR model \footnote{In \citep{daganzo2006variational}, the inhomogeneous LWR model is solved in the space of cumulative number of vehicles as a calculus of variations problem, and the existence of its solution is proved for road links with point bottlenecks. However, the wave solutions of the Riemann problem are not explicitly discussed.}.
In the first type, the inhomogeneous LWR model can be analyzed as a non-strictly hyperbolic conservation law \citep{isaacson1992resonance,lin1995resonant,jin2003inhLWR} or as a hyperbolic conservation law with a discontinuous flux function \citep{gimse1990discontinuous,gimse1993discontinuous,klingenberg1995discontinuous,diehl1995discontinuous,diehl1996discontinuous,diehl1996discontinuous2,diehl1996discontinuous3,zhang2003flux,burger2005discontinuous,burger2008inhomogeneous}, and various numerical methods can be used \citep{bale2002flux,zhang2005inhomogeneous,zhang2005fluxes,zhang2006inhomogeneous,herty2007discontinuous}.
In the second type, the self-similar waves of the Riemann solutions are separated into links 1 and 2 by introducing a stationary state for each link, and the wave on each link is determined by a new Riemann problem of the corresponding homogeneous LWR model \citep{seguin2003discontinuous,garavello2007discontinuous}. Here the stationary states are subject to admissible conditions as well as certain entropy conditions. This solution framework was first proposed for solving Riemann problems at general junctions with more than one upstream and downstream link \citep{holden1995unidirection,coclite2005network}. In \citep{seguin2003discontinuous}, the method was introduced for solving the inhomogeneous LWR model, and the stationary states are solved for a specific example. In \citep{garavello2007discontinuous}, a more general approach was proposed for solving the stationary states with a singular map method.
However, all these existing methods solve the Riemann problem in flux-density space: the first type of  method is tedious due to the need to analyze kinematic waves on both links at the same time, and the second type of
 method fails to present the entropy condition in a physically meaningful way. In addition, all existing methods do not account for interior states in stationary shock waves \citep{vanleer1984upwind,bultelle1998shock} and cannot be easily extended for studying traffic dynamics in a road network \citep{jin2003dissertation}. Note that, in this paper, we do not intend to study numerical solution methods for solving the inhomogeneous LWR model.

In this paper, we present a new framework for analyzing the inhomogeneous LWR model. We also adopt the method of wave separation by \citep{holden1995unidirection}, but introduce a stationary state and an interior state for each branch. Here stationary states are the self-similar states at the boundary, and interior states do not take any space in the continuous solution and only show up in the numerical solutions as observed in \citep{vanleer1984upwind,bultelle1998shock}. Rather than using the fundamental diagram, we introduce a so-called supply-demand diagram and discuss the problem in supply-demand space. After deriving admissible solutions for upstream and downstream stationary and interior states in supply-demand space, we introduce an entropy condition based on the discrete supply-demand method \citep{daganzo1995ctm,lebacque1996godunov}. We then prove that stationary states exist and are unique for given upstream demand and downstream supply, and interior states exist but may not be unique. Further we compare the Riemann solutions obtained by the new method with those in the literature for both the homogeneous and inhomogeneous LWR models. We also apply the new framework for analyzing asymptotic stationary states on an inhomogeneous ring road and demonstrate the existence of interior states with numerical examples.

The rest of the paper is organized as follows. In Section 2, we review the definitions of the supply and demand functions and the discrete supply-demand method for computing boundary fluxes. In Section 3, we introduce the supply-demand diagrams and the structure of the solutions to the Riemann problem of the inhomogeneous LWR model in supply-demand space. In section 4, we derive the admissible conditions for stationary and interior states in  supply-demand space and an entropy condition consistent with the local supply-demand method in interior states. In Section 5, we solve the Riemann problem for both the homogeneous and inhomogeneous LWR models and present a graphical solution scheme. In Section 6, we analyze asymptotic stationary states on an inhomogeneous ring road and demonstrate the existence of interior states with numerical solutions. In Section 7, we conclude our study with a discussion of future directions.

\section{Review of the supply-demand functions and methods}

\subsection{Review of Engquist-Osher functions and the Godunov method for convex conservation laws}
For the original LWR model \refe{lwr}, assuming that $k=\r_c-\r$, we obtain a hyperbolic conservation law in $k\in [-\r_c, \r_j-\r_c]$ as follows
\bqn
\pd{k}{t}+\pd{f(k)}x&=&0, \label{lwr:convex}
\eqn
where $f(k)=C-Q(\r_c-k)$ is convex when $Q(\r)$ is concave, since $\pdd{f(k)}k=-\pdd{Q(\r)}{\r}\geq 0$. Here $f(0)=0$. Moreover, if $q=Q(\r)$ is the Greenshields fundamental diagram \citep{greenshields1935capacity}, then \refe{lwr:convex} is Burgers' equation.

For the nonlinear equation \refe{lwr:convex}, we usually have to resort to numerical solutions for general initial and boundary conditions. After dividing the time duration into a number of time intervals of $\dt$ and splitting the road link into a number of cells of width $\dx$, the finite difference equation in conservation form can be written as follows \citep{colella2000numerical,leveque2002fvm}:
\bqn
k_i^{j+1}=k_i^j-\frac{\dt}{\dx} \left(f_{i+1/2}^{j*}-f_{i-1/2}^{j*} \right), \label{finitedifference}
\eqn
where $k_i^j$ is the average value of $k(x,t)$ in cell $i$ between $x_{i}-\frac{\dx}2$ and $x_{i}+\frac{\dx}2$ at time step $j$, and $f_{i+1/2}^{j*}$ is the flux through the boundary $x_1+\frac{\dx}2$ between time steps $j$ and $j+1$. Here $\dx$ and $\dt$ have to satisfy the so-called CFL condition \citep{courant1928CFL}.

For a hyperbolic conservation law \refe{lwr:convex}, the following functions were first introduced by \citep{engquist1980difference,engquist1980calculation,engquist1981difference,osher1982upwind}
\bqn
g(k)&=&f(\max\{k,0\})=\cas{{ll}f(k),&\m{if } k\geq 0\\0,&\m{if }k\leq 0},\nonumber \\&=&\int_0^k \chi(s) f'(s) ds=\int_0^k \max\{f'(s),0\} ds,\\
h(k)&=&f(\min\{k,0\})=\cas{{ll}f(k),&\m{if } k\leq 0\\0,&\m{if }k\geq 0},\nonumber\\&=&\int_0^k (1- \chi(s)) f'(s) ds=\int_0^k \min\{f'(s),0\} ds,
\eqn
where $\chi(k)$ equals 1 iff $f'(k)\geq 0$ and equals 0 otherwise. Note that $f(k)=g(k)+h(k)$ and $f(k)=\max\{g(k),h(k)\}$.
Therefore, we can rewrite \refe{lwr:convex} in the following form:
		\bqn
		k_t+g(k)_x+h(k)_x&=&0,\label{lwr:fluxsplit}\\
		k_t+\left[ \max\{g(k),h(k)\}\right]_x&=&0.
		\eqn

Further, based on these definitions, the following E-O flux is introduced \citep{engquist1980calculation,osher1984difference}
\bqn
f_{i-1/2}^{j*}&=&g(k_{i-1}^j)+h(k_{i}^j)=f(k_{i-1}^j)+\int_{k_{i-1}}^{k_i} \min\{f'(s), 0\}ds\nonumber \\&=&f(k_{i}^j)-\int_{k_{i-1}}^{k_i} \max\{f'(s), 0\}ds \nonumber\\
&=&\frac 12\left[f(k_{i-1}^j)+f(k_{i}^j)-\int_{k_{i-1}}^{k_i} |f'(s)|ds \right].
\eqn
That is, \refe{finitedifference} can be written as
\bqn
k_i^{j+1}=k_i^j-\frac{\dt}{\dx} \left( h(k_{i+1}^j)-h(k_{i}^j)+g(k_{i}^j)-g(k_{i-1}^j) \right).
\eqn
This can be considered as upwind method for \refe{lwr:fluxsplit}, since $g(k)$ is non-decreasing and $h(k)$ non-increasing.

In \citep{vanleer1984upwind}, the Godunov flux \citep{godunov1959} for Burgers' equation was written as
\bqn
f_{i-1/2}^{j*}&=&\max\{g(k_{i-1}^j),h(k_{i}^j)\}.  \label{sdflux}
\eqn
In \citep{osher1984difference}, a new formulation of the Godunov flux \citep{godunov1959} was introduced as
\bqn
f_{i-1/2}^{j*}&=&\cas{{ll} \min_{k_{i-1}^j\leq k\leq k_i^j} f(k), & k_{i-1}^j\leq k_i^j,\\
\max_{k_{i-1}^j\geq k\geq k_i^j} f(k), & k_{i-1}^j> k_i^j.} \label{osherflux}
\eqn
For convex $f(k)=\max\{g(k),h(k)\}$, since $g(k)$ and $h(k)$ are monotonically increasing and decreasing respectively, this is equivalent to
\bqs
f_{i-1/2}^{j*}&=&\cas{{ll} \max\{ \min_{k_{i-1}^j\leq k\leq k_i^j} g(k),\min_{k_{i-1}^j\leq k\leq k_i^j} h(k) \}, & k_{i-1}^j\leq k_i^j\\
\max\{\max_{k_{i-1}^j\geq k\geq k_i^j} g(k), \max_{k_{i-1}^j\geq k\geq k_i^j} h(k)\}, & k_{i-1}^j> k_i^j}\\&=&\max\{g(k_{i-1}^j),h(k_{i}^j)\}.
\eqs
That is, \refe{osherflux} is equivalent to \refe{sdflux}. However, it has been shown that \refe{osherflux} can also be applied to non-convex $f(k)$.

\subsection{Review of supply and demand functions and Godunov methods for the LWR model}
For the LWR model \refe{lwr}, we define the following functions
\bqn
D(\r)&=&Q(\min\{\r,\r_c\})=\cas{{ll}Q(\r),&\m{if } \r\leq \r_c\\C,&\m{if }\r\geq \r_c},\nonumber \\&=&\int_0^\r \chi(s) Q'(s) ds=\int_0^\r \max\{Q'(s),0\} ds \label{def:d}\\
S(\r)&=&Q(\max\{\r,\r_c\})=\cas{{ll}Q(\r),&\m{if } \r\geq \r_c\\C,&\m{if }\r\leq \r_c},\nonumber\\&=&C+\int_0^\r (1- \chi(s)) Q'(s) ds=C+\int_0^\r \min\{Q'(s),0\} ds, \label{def:s}
\eqn
where $\chi(\r)$ equals 1 iff $Q'(\r)\geq 0$ and equals 0 otherwise. It is straightforward to show that $D(\r)=C-g(k)$ and $S(\r)=C-h(k)$. Therefore, the Godunov method for \refe{lwr} is equivalent to
\bqn
\r_i^{j+1}=\r_i^j-\frac{\dt}{\dx} \left(q_{i+1/2}^{j*}-q_{i-1/2}^{j*} \right), \label{lwrfd}
\eqn
where the boundary flux can be written as \citep{vanleer1984upwind}
\bqn
q_{i-1/2}^{j*}&=&\min\{D(\r_{i-1}^j),S(\r_{i}^j)\}, \label{sdmethod}
\eqn
or as \citep{osher1984difference}
\bqn
q_{i-1/2}^{j*}&=&\cas{{ll} \min_{\r_{i-1}^j\leq\r\leq \r_i^j} Q(\r), &
\r_{i-1}^j<\r_i^j\\\max_{\r_i^j\leq \r\leq \r_{i-1}^j} Q(\r), & \r_i^j<\r_{i-1}^j }. \label{oshermethod}
\eqn
From \refe{oshermethod}, we can see that \refe{sdmethod} is valid as long as the fundamental diagram $Q(\r)$ is unimodal and may not be concave.

In the transportation literature, \citep{bui1992entropy} first applied \refe{lwrfd} and \refe{oshermethod} for solving the LWR model.
In \citep{daganzo1995ctm}, a new finite difference form was proposed for the LWR model with a triangular or trapezoidal fundamental diagram
\bqn
Q(\r)=\min \{v_f\r, v_c(\r_j-\r), Q_{max} \}
\eqn
where $v_c=v_f \frac{\r_c}{\r_j-\r_c}$ is the absolute value of the shock wave speed in congested traffic.
In the so-called cell transmission model (CTM), the space-time domain was discretized with a CFL number of 1; i.e., $\dx=v_f \dt$, and the boundary flux in \refe{lwrfd} was written as
\bqn
q_{i-1/2}^{j*}&=&\min\{\bar D(\r_{i-1}^j),\bar S(\r_{i}^j)\}, \label{ctmodel}
\eqn
where $\bar D(\r_i^j)\dt=\min \{Q_{max}\dt, n_i^j\}$ and $\bar S(\r_i^j)\dt=\min \{Q_{max}\dt, \frac{v_c}{v_f} (N_{max}-n_i^j)\} $ are defined as ``the maximum flows that can be sent and received by cell $i$ in the interval between time steps $j$ and $j+1$", $n_i^j=\r_i^j \dx=\r_i^j v_f \dt$ is the number of vehicles in cell $i$ at time step $j$, and $N_{max}=\r_j \dx=\r_j v_f\dt$ is the maximum number of vehicles in cell $i$.
Hence the physical meaning of \refe{ctmodel} is that the boundary flux is the minimum of the upstream sending flux and the downstream receiving flux. It can be shown that $\bar D(\r)=D(\r)$ and $\bar S(\r)=S(\r)$. Thus, \refe{ctmodel} is equivalent to \refe{sdmethod} for a CFL number of 1 and triangular or trapezoidal fundamental diagrams.

Following \citep{lebacque1996godunov}, we refer to $D(\r)$ in \refe{def:d} and $S(\r)$ in \refe{def:s} as the demand and supply functions respectively and call \refe{sdmethod} the discrete supply-demand method for computing fluxes.
The physical interpretations of demand and supply functions and the supply-demand method have formed the basis for extending the supply-demand method for computing fluxes through various network junctions \citep{daganzo1995ctm,lebacque1996godunov,jin2003dissertation}.
For inhomogeneous roads, the extension is straightforward as follows \citep{daganzo1995ctm,lebacque1996godunov}: the demand and supply functions in \refe{def:d} and \refe{def:s} are location-dependent, $D(x,\r(x,t))=Q(x,\min\{\r(x,t),\r_c(x)\})$ and $S(x,\r(x,t))=Q(x,\max\{\r(x,t),\r_c(x)\})$, and the flux is still computed by the supply-demand method \refe{sdmethod}, $q_{i-1/2}^{j*}=\min\{D_{i-1}(\r_{i-1}^j),S_i(\r_{i}^j)\}$, where $D_{i-1}(\r)$ is the demand function in cell $i-1$, and $S_i(\r)$ is the supply function in cell $i$. It has been shown that the flux by the extended supply-demand method is still the Godunov flux \citep{jin2003inhLWR,zhang2003flux}.

In this study, based on the Godunov finite difference equation in \refe{lwrfd} and the supply-demand method in \refe{sdmethod} for computing boundary fluxes, we attempt to construct the convergent solution of \refe{inhlwr} with discontinuous flux functions \refe{discontinuousflux} and initial conditions \refe{inicon}.
In \citep{daganzo1995difference}, the convergence, truncation error, and capture of shock waves were directly derived from the corresponding finite difference equation of the homogeneous LWR model.
Here, we attempt to present a new framework and find the solutions to the Riemann problem of the inhomogeneous LWR model at a linear boundary.

\section{Supply-demand diagrams and the structure of Riemann solutions}
In the literature, the inhomogeneous LWR model has been analyzed in flux-density space. Since a traffic state in flux-density space $(\r, q)$ can also be represented in supply-demand space as $U=(D,S)$, in this study we will analyze the inhomogeneous LWR model in supply-demand space.
Furthermore, we present the structure of solutions to the Riemann problem of the inhomogeneous LWR model in supply-demand space.

\subsection{Supply-demand diagrams}
Corresponding to the fundamental diagram in flux-density space, a supply-demand diagram can be introduced in supply-demand space.
In \reff{fig:fd-ds}(b), we draw a supply-demand diagram for the two fundamental diagrams in \reff{fig:fd-ds}(a).
On the dashed branch of the supply-demand diagram, traffic is under-critical (UC) and $U=(D,C)$ with $D\leq C$; on the solid branch, traffic is over-critical (OC) and $U=(C,S)$ with $S\leq C$.
Compared with the fundamental diagram of a road section, the supply-demand diagram only considers capacity $C$ and congestion level of traffic flow, but not other detailed characteristics such as  critical density, jam density, or relationship between density and flux.
That is, different fundamental diagrams can have the same supply-demand diagram, as long as they have the same capacity and are unimodal, and their critical densities, jam densities, or shapes are not relevant.
However, there is always a one-to-one mapping between a given supply-demand diagram and its corresponding fundamental diagram.
That is, there exists a one-to-one mapping between $(\rho,q)$ and $(D,S)$.

\bfg\bc $\ba{c@{\hspace{0.3in}}c}
\includegraphics[height=1.3in]{sta20061026figure.2} &
\includegraphics[height=1.3in]{sta20061026figure.1} \\
\multicolumn{1}{c}{\mbox{\bf (a)}} &
    \multicolumn{1}{c}{\mbox{\bf (b)}}
\ea$ \ec \caption{The fundamental diagram and its supply-demand diagram}\label{fig:fd-ds} \efg

For the demand and supply functions in \refe{def:d} and \refe{def:s}, we can see that $D$ is non-decreasing with $\r$ and $S$ non-increasing.
Thus $D\leq C$, and $S\leq C$, and $\max \{D,S\}=C$. In addition, $D=S=C$ iff traffic is critical; $D<S=C$ iff traffic is strictly under-critical (SUC); $S<D=C$ iff traffic is strictly over-critical (SOC). Therefore, the state $U=(D,S)$ is under-critical (UC), iff $S=C$, or equivalently $D\leq S$. The state $U=(D,S)$ is over-critical (OC), iff $D=C$, or equivalently $S\leq D$.

From a state on the supply-demand diagram, we can obtain the corresponding flux $q(U)=\min\{D,S\}$ and capacity $C=\max \{D,S\}$. However, we cannot tell the density from the supply-demand diagram, and the fundamental diagram is still needed to compute the density $\rho$ from $(D,S)$. That is, $\rho$ can be written as a function of $(D,S)$ as
\bqs
\rho(D,S)&=&\cas{{ll} D^{-1}(D), &D\leq S,\\ S^{-1}(S), &S<D,}
\eqs
since $D(\r)$ and $S(\r)$ are invertible when the traffic is UC and OC respectively.
Note that $\rho$ is not a function of $q$.
If we introduce the supply-demand ratio $\gamma=D/S$, then $q(D,S)=min\{\gamma,1/\gamma\}\cdot C$, and
\bqs
\r(D,S)&=&R(\gamma)\equiv \cas{{ll} D^{-1}(C\gamma), &\gamma\leq 1, \\S^{-1}(\frac{C}{\gamma}), &\gamma >1,}
\eqs
 where $R(\gamma)$ is an increasing function in $\gamma\in[0,\infty]$. Here $R(0)=0$, $R(1)=\r_c$, and $R(\infty)=\r_j$. In this sense, $\r=R(\gamma)$ can be considered as the inverse flux-density relationship. Similarly, $v=V(\r)=V(R(\gamma))$ is a non-increasing function in $\gamma$, $V(0)=v_f$, and $V(\infty)=0$.

\subsection{The structure of solutions to the Riemann problem}
In supply-demand space, initial conditions \refe{inicon} are equivalent to
\bqn
U(x,t=0)&=&\cas{{ll} U_1=(D_1,S_1), &x<0, \\ U_2=(D_2,S_2), &x>0.} \label{dsriemann}
\eqn
The Riemann problem at a linear boundary is then equivalent to the Riemann problem for \refe{inhlwr} with initial conditions \refe{dsriemann}.

Unlike existing studies of hyperbolic conservation laws with discontinuous flux functions, in which solutions to the Riemann problem have been constructed in $\r-q$ space, within the framework of wave separation developed in \citep{holden1995unidirection}, we construct the solutions to the Riemann problem for \refe{inhlwr} with initial conditions \refe{dsriemann} in supply-demand space.

In solutions to the Riemann problem for \refe{inhlwr} with initial conditions \refe{dsriemann}, a shock wave or a rarefaction wave could initiate on a link from the linear boundary $x=0$, and traffic states on both links become asymptotically stationary after a long time.
We denote the stationary state on the upstream link $1$ by $U_1^-$ and the stationary state on the downstream link $2$ by $U_2^+$.
At the boundary, there can also exist interior states \citep{vanleer1984upwind,bultelle1998shock}, which take infinitesimal space. We denote the interior states on links 1 and 2 by $U_1(0^-,t)$ and $U_2(0^+,t)$ respectively. The structure of Riemann solutions on upstream and downstream links are shown in \reff{fig:riemannstructure}.

\bfg\bc $\ba{c@{\hspace{0.3in}}c}
\includegraphics[width=2in]{sd20080619figure.17} &
\includegraphics[width=2in]{sd20080619figure.18} \\
\multicolumn{1}{c}{\mbox{\bf (a)}} &
    \multicolumn{1}{c}{\mbox{\bf (b)}}
\ea$ \ec \caption{Structure of Riemann solutions: (a) Upstream link 1; (b) Downstream link 2}\label{fig:riemannstructure} \efg

Then the kinematic wave on the upstream link $1$ is the solution of the corresponding homogeneous LWR model with initial left and right conditions of $U_1$ and $U_1^-$, respectively. Similarly, the kinematic wave on the downstream link $2$ is the solution of the corresponding LWR model with initial left and right conditions of $U_2^+$ and $U_2$, respectively.
Since the stationary and interior states are constant for $t>0$, states on both links 1 and 2 are self-similar \citep{smoller1983shock}. That is, if stationary states exist and are unique, we have unique self-similar solutions for the Riemann problem of \refe{inhlwr}.
In the following sections, we first derive necessary conditions for both stationary and interior states and then solve them in supply-demand space.

\section{Necessary conditions for the existence of stationary and interior states}
We denote $q_{1\to 2}(0,t)$ as the flux from link $1$ to link $2$. We first observe that the fluxes are determined by the stationary states: the asymptotic out-flux of link 1 is $q_1(0^-,t)=q(U_1^-)$, and the asymptotic in-flux of link 2 is $q_2(0^+,t)=q(U_2^+)$.
Furthermore, from the conservation of traffic at the linear boundary, we have
\bqn\label{trafficconservation}
q_{1\to 2}(0,t)&=&q_1(0^-,t)=q_2(0^+,t)=q(U_1^-)=q(U_2^+).
\eqn

\subsection{The admissible conditions for stationary states}
 As observed in \citep{holden1995unidirection,coclite2005network}, the speed of the kinematic wave on an upstream link cannot be positive, and that on a downstream link cannot be negative. We have the following theorem.

\bfg\bc $\ba{c@{\hspace{0.3in}}c}
\includegraphics[height=1.3in]{sd20080619figure.1} &
\includegraphics[height=1.3in]{sd20080619figure.2} \\
\multicolumn{1}{c}{\mbox{\bf (a)}} &
    \multicolumn{1}{c}{\mbox{\bf (b)}}\\
\includegraphics[height=1.3in]{sd20080619figure.3} &
\includegraphics[height=1.3in]{sd20080619figure.4} \\
\multicolumn{1}{c}{\mbox{\bf (c)}} &
    \multicolumn{1}{c}{\mbox{\bf (d)}}	
\ea$ \ec \caption{Admissible stationary states for the upstream link $1$}\label{fig:ssupstream} \efg

\bfg\bc $\ba{c@{\hspace{0.3in}}c}
\includegraphics[height=1.3in]{sd20080619figure.5} &
\includegraphics[height=1.3in]{sd20080619figure.6} \\
\multicolumn{1}{c}{\mbox{\bf (a)}} &
    \multicolumn{1}{c}{\mbox{\bf (b)}}\\
\includegraphics[height=1.3in]{sd20080619figure.7} &
\includegraphics[height=1.3in]{sd20080619figure.8} \\
\multicolumn{1}{c}{\mbox{\bf (c)}} &
    \multicolumn{1}{c}{\mbox{\bf (d)}}	
\ea$ \ec \caption{Admissible stationary states for the downstream link $2$}\label{fig:ssdownstream} \efg

\begin{theorem} [Admissible stationary states] \label{admissibless} For initial conditions in \refe{dsriemann},  stationary states are admissible if and only if
\bqn
U_1^-&=&(D_1,C_1) \m{ or } (C_1,  S_1^-), \label{upstreamss}
\eqn
where $S_1^-<D_1$, and
\bqn
U_2^+&=&(C_2,S_2) \m{ or } (D_2^+, C_2), \label{downstreamss}
\eqn
where $D_2^+<S_2$. The results are illustrated in Figures \ref{fig:ssupstream} and \ref{fig:ssdownstream}.
\end{theorem}
{\em Proof.}
When the initial state on link $1$ is strictly under-critical; i.e., when $D_1<S_1=C_1$, the admissible stationary state $U_1^-$ is the same as $U_1=(D_1,C_1)$ or strictly over-critical with $U_1^-=(C_1, S_1^-)$, where $S_1^-<D_1$. In this case, the Riemann problem for the LWR model on the upstream link $1$ with upstream and downstream initial states $U_1=(D_1, C_1)$ and $U_1^-$ has the following possible solutions: there is no wave when $U_1^-=U_1$; there is a backward traveling shock wave when $U_1^-=(C_1, S_1^-)$. In addition, we can verify that any stationary states not equal to $U_1$ and with $S_1^->D_1$ will lead to forward traveling shock waves or rarefaction waves. Note that when $U_1^-=(C_1, D_1)$, the Riemann problem is solved by a zero shock wave, but $U_1^-$ cannot be the stationary state by definition.
When the initial state on link $1$ is over-critical; i.e., when $S_1\leq D_1 =C_1$, the admissible stationary state $U_1^-$ is over-critical with $U_1^-=(C_1, S_1^-)$, where $S_1^-\leq C_1$. In this case, the Riemann problem for the LWR model on the upstream link $1$ with upstream and downstream initial states $U_1=(C_1, S_1)$ and $U_1^-=(C_1, S_1^-)$ has the following possible solutions: there is no wave when $S_1^-=S_1$; there is a backward traveling shock wave when $S_1>S_1^-$; and there is a backward traveling rarefaction wave when $S_1<S_1^-$. Therefore, the stationary state is indeed admissible. In addition, we can verify that any strictly under-critical stationary states $U_1^-$ will lead to forward traveling rarefaction waves and are not admissible.

When the initial state on link $2$ is under-critical; i.e., when $D_2\leq S_2=C_2$, the admissible stationary state $U_2^+$ is under-critical with $U_2^+=(D_2^+, C_2)$, where $D_2^+\leq C_2$. In this case, the Riemann problem for the LWR model on the downstream link $2$ with upstream and downstream initial states $U_2^+=(D_2^+, C_2)$ and $U_2=(D_2,C_2)$ has the following possible solutions: there is no wave when $D_2^+=D_2$; there is a forward traveling shock wave when $D_2^+<D_2$; there is a forward traveling rarefaction wave when $D_2^+>D_2$. Therefore, the stationary state is indeed admissible. In addition, we can verify that any strictly over-critical stationary states $U_2^+$ will lead to backward traveling rarefaction waves and are not admissible. When the initial state on link $2$ is strictly over-critical; i.e., when $S_2< D_2=C_2$, the admissible stationary state $U_2^+$ is the same as $U_2=(C_2, S_2)$ or strictly under-critical with $U_2^+=(D_2^+, C_2)$, where $D_2^+ < S_2$. In this case, the Riemann problem for the LWR model on the downstream link $2$ with upstream and downstream initial states $U_2^+$ and $U_2=(D_2,C_2)$ has the following possible solutions: there is no wave when $U_2^+=U_2$; there is a forward traveling shock wave when $U_2^+=(D_2^+, C_2)$. In addition, we can verify that any stationary states not equal to $U_2$ and with $D_2^+>S_2$ will lead to backward traveling shock waves or rarefaction waves. Note that when $U_2^+=(S_2, C_2)$, the Riemann problem is solved by a zero shock wave, but $U_2^+$ cannot be the stationary state by definition.
\eop

{\em Remark 1.} Note that $U_1^-=U_1$ and $U_2^+=U_2$ are always admissible. In this case, the stationary states are the same as the corresponding initial states, and there are no waves.

{\em Remark 2.} From the proof of Theorem \ref{admissibless}, we can see that the types and traveling directions of waves on a homogeneous road can be solely determined by upstream demand and downstream supply and are not related to the shape of fundamental diagrams. This is why we are able to discuss the Riemann problem of the inhomogeneous LWR model in supply-demand space. Note that, however, the wave speeds are related to the details of the flux-density relation $Q(\r)$.

{\em Remark 3.}
Then the out-flux $q_1(0^-,t)=q(U_1^-)=\min\{D_1^-,S_1^-\}\leq D_1$ and the in-flux $q_2(0^+,t)=q(U_2^+)=\min\{D_2^+,S_2^+\}\leq S_2$.
That is, $D_1$ is the maximum sending flux and $S_2$ is the maximum receiving flux in the sense of \citep{daganzo1994ctm,daganzo1995ctm}.
Furthermore, $q_1(0^-,t)=D_1$, iff $U_1^-=(D_1,C_1)$, and iff $U_1^-$ is UC; $q_1(0^-,t)<D_1$ iff $U_1^-=(C_1,S_1^-)$ with $S_1^-<D_1$, and iff $U_1^-$ is SOC. Similarly, $q_2(0^+,t)=S_2$, iff $U_2^+=(C_2,S_2)$, and iff $U_2^+$ is OC; $q_2(0^+,t)<S_2$ iff $U_2^+=(D_2^+,C_2)$ with $D_2^+<S_2$, and iff $U_2^+$ is SUC.

\subsection{The admissible conditions for interior states}
The Riemann problem on link 1 with left and right initial conditions of $U_1^-$ and  $U_1(0^-,t)$ cannot have negative waves. Otherwise, $U_1(0^-,t)$ will propagate upstream and violates the condition that it only exists at the boundary, but not anywhere upstream. Similarly, the Riemann problem on link 2 with left and right initial conditions of $U_2(0^+,t)$ and $U_2^+$ cannot have positive waves. Therefore, interior states $U_1(0^-,t)$ and $U_2(0^+,t)$ should satisfy the following admissible conditions.

\bfg\bc $\ba{c@{\hspace{0.3in}}c}
\includegraphics[height=1.3in]{sd20080619figure.9} &
\includegraphics[height=1.3in]{sd20080619figure.10} \\
\multicolumn{1}{c}{\mbox{\bf (a)}} &
    \multicolumn{1}{c}{\mbox{\bf (b)}}\\
\includegraphics[height=1.3in]{sd20080619figure.11} &
\includegraphics[height=1.3in]{sd20080619figure.12} \\
\multicolumn{1}{c}{\mbox{\bf (c)}} &
    \multicolumn{1}{c}{\mbox{\bf (d)}}	
\ea$ \ec \caption{Admissible interior states for upstream link $1$}\label{fig:isupstream} \efg

\bfg\bc $\ba{c@{\hspace{0.3in}}c}
\includegraphics[height=1.3in]{sd20080619figure.13} &
\includegraphics[height=1.3in]{sd20080619figure.14} \\
\multicolumn{1}{c}{\mbox{\bf (a)}} &
    \multicolumn{1}{c}{\mbox{\bf (b)}}\\
\includegraphics[height=1.3in]{sd20080619figure.15} &
\includegraphics[height=1.3in]{sd20080619figure.16} \\
\multicolumn{1}{c}{\mbox{\bf (c)}} &
    \multicolumn{1}{c}{\mbox{\bf (d)}}	
\ea$ \ec \caption{Admissible interior states for downstream link $2$}\label{fig:isdownstream} \efg

\begin{theorem}[Admissible interior states]\label{intersta}
For asymptotic stationary states $U_1^-$ and $U_2^+$, interior states $U_1(0^-,t)$ and $U_2(0^+,t)$ are admissible if and only if
\bqn
U_1(0^-,t)&=&\cas{{ll}(C_1,S_1^-)=U_1^-, &\m{when }S_1^-< D_1^-=C_1 \\(D_1(0^-,t),  S_1(0^-,t)), & \m{when }D_1^-\leq S_1^-=C_1} \label{upstreamis}
\eqn
where $S_1(0^-,t) \geq D_1^-$, and
\bqn
U_2(0^+,t)&=&\cas{{ll}(D_2^+,C_2)=U_2^+,&\m{when }D_2^+<S_2^+=C_2 \\ (D_2(0^+,t), S_2(0^+,t)), &\m{when }S_2^+\leq D_2^+=C_2 } \label{downstreamis}
\eqn
where $D_2(0^+,t)\geq S_2^+$. The results are illustrated in Figures \ref{fig:isupstream} and \ref{fig:isdownstream}.
\end{theorem}
{\em Proof}. The results can be verified with the observation that the Riemann solutions for the homogeneous LWR model of the upstream links cannot have negative waves, and those of the downstream links cannot have positive waves. \eop

{\em Remark 1.} Note that $U_1(0^-,t)=U_1^-$ and $U_2(0^+,t)=U_2^+$ are always admissible. In this case, the interior states are the same as stationary states, and it is equivalent to saying that there are no interior states.

\subsection{An entropy condition for the local supply-demand method}
In addition to traffic conservation and admissible conditions, we introduce an entropy condition such that the boundary flux is always consistent with that by the supply-demand method \refe{sdmethod} for the local interior states. At the boundary at $x=0$, the immediate upstream state is $U_1(0^-,t)$, and the immediate downstream state $U_2(0^+,t)$. That is, the interior states have to satisfy the following entropy condition
\bqn
q_{1\to 2}(0,t) &=& \min\{D_1(0^-,t),S_2(0^+,t) \}.\label{sdentropy}
\eqn

Note that the entropy condition \refe{sdentropy} is also equivalent to the following localized optimization problem
\bqs
\max \{q_{1\to 2}(0,t)=q_1(0^-,t)=q_2(0^+,t)\}
\eqs
subject to
\bqs
q_1(0^-,t) &\leq & D_1(0^-,t),\\
q_2(0^+,t) &\leq & S_2(0^+,t).
\eqs
That is, the stationary and interior states are solutions of the optimization problem in the domain defined by the traffic conservation condition \refe{trafficconservation}, the admissible conditions for stationary states (\ref{upstreamss}-\ref{downstreamss}), and the admissible conditions for interior states (\ref{upstreamis}-\ref{downstreamis}).
Optimization formulations of entropy conditions were also adopted in \citep{holden1995unidirection,coclite2005network}, but in terms of stationary states in flux-density space rather than in supply-demand space as we have done here.

In all of the necessary conditions above, we can see that the stationary and interior states are independent of the upstream supply, $S_1$, and the downstream demand, $D_2$. That is, the same upstream demand and downstream supply will yield the same solutions of stationary and interior states: when the upstream traffic is congested, its congestion level is not relevant to the stationary and interior states or the boundary flux; when the downstream traffic is free flow, its density is not relevant to the stationary and interior states or the boundary flux.
Note that, however, the types and speeds of waves on both links can be related to $S_1$ as shown in \reff{fig:ssupstream}(d) and $D_2$ as shown in \reff{fig:ssdownstream}(d).

\section{Solutions to the Riemann problem for the inhomogeneous LWR model}
For the Riemann problem of \refe{inhlwr} with initial conditions \refe{dsriemann}, we first solve for the stationary and interior states that satisfy the traffic conservation condition \refe{trafficconservation}, the admissible conditions for stationary states (\ref{upstreamss}-\ref{downstreamss}), the admissible conditions for interior states (\ref{upstreamis}-\ref{downstreamis}), and the entropy condition \refe{sdentropy}.
Then we compare them with existing solutions for both the homogeneous and inhomogeneous LWR models.

\begin{lemma} \label{globalflux} In the Riemann solutions, the boundary flux satisfies
\bqn
q_{1\to 2}(0,t)=q(U_1^-)=q(U_2^+) =  \min\{D_1, S_2\}. \label{flowsta}
\eqn
\end{lemma}
{\em Proof.} From Theorem \ref{admissibless} and traffic conservation \refe{trafficconservation}, we have that
\bqs
q_{1\to 2}(0,t)=q(U_1^-)=q(U_2^+) \leq  \min\{D_1, S_2\}.
\eqs
We demonstrate that it is not possible for $q_{1\to 2}(0,t)<\min\{D_1, S_2\}$. Otherwise, we assume that $q_{1\to 2}(0,t)=q_0<\min\{D_1, S_2\}\leq \min\{C_1, C_2\}$. Since $q(U_1^-)=q_0<D_1$, from \refe{upstreamss} we have that $U_1^-=(C_1,q_0)$. Further from \refe{upstreamis} we have that $U_1(0^-,t)=U_1^-=(C_1,q_0)$. Hence $D_1(0^-,t)=C_1$.
Similarly, since $q(U_2^+)=q_0<S_2$, from \refe{downstreamss} we have that $U_2^+=(q_0, C_2)$. Further from \refe{downstreamis} we have that $U_2(0^+,t)=U_2^+=(q_0,C_2)$. Hence $S_2(0^+,t)=C_2$.
Then from \refe{sdentropy} we have $q_{1\to 2}(0,t)=\min\{C_1, C_2\}$, which contradicts the assumption that $q_{1\to 2}(0,t)<\min\{C_1, C_2\}$.
Therefore we have \refe{flowsta}.

That is, the local optimal solution in \refe{sdentropy} leads to a global optimal flux at a linear boundary, which satisfies the following optimization problem
\bqs
\max \{q_{1\to 2}(0,t)\}
\eqs
subject to
\bqs
q_{1\to 2}(0,t) &\leq & D_1,\\
q_{1\to 2}(0,t) &\leq & S_2.
\eqs

Therefore, in the Godunov finite difference equation \refe{lwrfd}, the boundary flux at the first time step by \refe{sdmethod} is the same as the asymptotic flux, regardless of the time step size. That is, the discrete flux in \refe{sdmethod} is consistent with the continuous flux in \refe{sdentropy}.
\eop

\begin{theorem} \label{linearboundary} The stationary states and interior states of the Riemann problem for \refe{inhlwr} with initial conditions \refe{dsriemann} are the following:
\ben
\item When $D_1 < S_2$, we have unique stationary and interior states: $U_1^-=U_1(0^-,t)=(D_1, C_1)$ and $U_2^+=U_2(0^+,t)=(D_1, C_2)$;
\item When $D_1 > S_2$, we have unique stationary and interior states: $U_1^-=U_1(0^-,t)=(C_1, S_2)$ and $U_2^+=U_2(0^+,t)=(C_2,S_2)$;
\item When $D_1 = S_2$, we have unique stationary states: $U_1^-=(D_1,C_1)$, and $U_2^+=(C_2,S_2)$; but interior states may not be unique: $U_1(0^-,t)=(D_1, C_1)$ and $D_2(0^+,t)\geq S_2$ and $S_2(0^+,t)\geq S_2$, or $U_2(0^+,t)=(C_2,S_2)$ and $S_1(0^-,t)\geq D_1$ and $D_1(0^-,t) \geq D_1$.
\een
\end{theorem}
{\em Proof.} When $D_1<S_2\leq C_2$, \refe{flowsta} leads to $q_{1\to 2}(0,t)=D_1$.
For link 1, since $q(U_1^-)=D_1$, from \refe{upstreamss} we have that $U_1^-=(D_1, C_1)$. Further from \refe{upstreamis} we have that $U_1(0^-,t)=(D_1(0^-,t), S_1(0^-,t))$ with $S_1(0^-,t)\geq D_1^-=D_1$.
For link 2, since $q(U_2^+)=D_1<S_2$, from \refe{downstreamss} we have that $U_2^+=(D_1, C_2)$. Further from \refe{downstreamis} we have that $U_2(0^+,t)=U_2^+=(D_1,C_2)$.
Then from \refe{sdentropy} we have $q_{1\to 2}(0,t)=\min\{D_1(0^-,t), C_2\}=D_1$. Therefore, $D_1(0^-,t)=D_1$ and $S_1(0^-,t)=C_1$.

When $S_2<D_1\leq C_1$, \refe{flowsta} leads to $q_{1\to 2}(0,t)=S_2$.
For link 1, since $q(U_1^-)=S_2<D_1$, from \refe{upstreamss} we have that $U_1^-=(C_1, S_2)$. Further from \refe{upstreamis} we have that $U_1(0^-,t)=U_1^-=(C_1,S_2)$.
For link 2, since $q(U_2^+)=S_2$, from \refe{downstreamss} we have that $U_2^+=(C_2, S_2)$. Further from \refe{downstreamis} we have that $U_2(0^+,t)=(D_2(0^+,t), S_2(0^+,t))$ with $D_2(0^+,t)\geq S_2$.
Then from \refe{sdentropy} we have $q_{1\to 2}(0,t)=\min\{C_1, S_2(0^+,t)\}=S_2$. Therefore, $S_2(0^+,t)=S_2$ and $D_2(0^+,t)=C_2$.

When $D_1=S_2$, \refe{flowsta} leads to $q_{1\to 2}(0,t)=D_1=S_2$.
For link 1, since $q(U_1^-)=D_1$, from \refe{upstreamss} we have that $U_1^-=(D_1, C_1)$. Further from \refe{upstreamis} we have that $U_1(0^-,t)=(D_1(0^-,t), S_1(0^-,t))$ with $S_1(0^-,t)\geq D_1^-=D_1$.
For link 2, since $q(U_2^+)=S_2$, from \refe{downstreamss} we have that $U_2^+=(C_2, S_2)$. Further from \refe{downstreamis} we have that $U_2(0^+,t)=(D_2(0^+,t), S_2(0^+,t))$ with $D_2(0^+,t)\geq S_2$.
Then from \refe{sdentropy} we have $q_{1\to 2}(0,t)=\min\{D_1(0^-,t), S_2(0^+,t)\}=D_1=S_2$.
If $D_1(0^-,t)=D_1=S_2$, then $S_1(0^-,t)=C_1$, and $S_2(0^+,t)\geq S_2$. In this case, the interior state $U_2(0^+,t)$ may not be unique with $D_2(0^+,t)\geq S_2$ and $S_2(0^+,t)\geq S_2$. Note that, when $S_2=C_2$, the interior state is $U_2(0^+,t)=U_2^+$.
If $S_2(0^+,t)=D_1=S_2$, then $D_2(0^+,t)=C_2$, and $D_1(0^-,t) \geq D_1$. In this case, the interior state $U_1(0^-,t)$ may not be unique with $S_1(0^-,t)\geq D_1$ and $D_1(0^-,t) \geq D_1$. Note that, when $D_1=C_1$, the interior state is $U_1(0^-,t)=U_1^-$.
If both the upstream link 1 and the downstream link 2 have the same fundamental diagram, this case corresponds to a stationary shock, and the interior state is the same as that in \citep{vanleer1984upwind}.
\eop

{\em Remark 1}. From the theorem we can see that the stationary states always exist and are unique for the same pair of $D_1$ and $S_2$. Thus, with given $U_1$ and $U_2$, we can find unique kinematic waves on both links 1 and 2. Therefore, in the new solution framework, the solutions for the Riemann problem of the inhomogeneous LWR model always exist and are unique, although we may have multiple interior states.

{\em Remark 2}. If the entropy condition \refe{sdentropy} is replaced by \refe{flowsta}, we still have the same solutions $U_1^-$ and $U_2^+$. That is, if we do not consider possible interior states as in \citep{seguin2003discontinuous,garavello2007discontinuous}, then traffic conservation \refe{trafficconservation}, admissible conditions for stationary states (\ref{upstreamss}-\ref{downstreamss}), and the entropy condition \refe{flowsta} will yield the same stationary state solutions. 
However, this simplified approach - which is what currently exists in the literature - does not yield the existence or properties of the interior states.

{\em Remark 3}. When $D_1=S_2$, we have the following interior states that are different from the stationary states at $x=0^-$ or $x=0^+$. The interior state at $x=0^-$ has to satisfy $q(U_1(0^-,t)) \geq D_1$, and the interior state at $x=0^+$ has to satisfy $q(U_2(0^+,t) \geq S_2$.

In addition, we have the following conclusions concerning possible stationary states at a linear boundary.
\begin{corollary}\label{linksta}
When both links 1 and 2 reach asymptotic stationary states, they share the same flux $q$, and possible stationary states are the following: both links are UC with link 1 at $(q,C_1)$ and link 2 at $(q, C_2)$ where $q=D_1<S_2$; both links are OC with link 1 at $(C_1, q)$ and link 2 at $(C_2, q)$ where $q=S_2<D_1$; link 1 is UC at $(q, C_1)$ and link 2 OC at $(C_2, q)$ where $q=D_1=S_2$. It is not possible that link 1 is SOC and link 2 SUC.
\end{corollary}
{\em Remark}. The stationary states are stable in the sense that, when they are given as initial states, we obtain the same stationary states following Theorem \ref{linearboundary}.

\subsection{The homogeneous LWR model}
\bfg \bc $\ba{c@{\hspace{0.3in}}c@{\hspace{0.3in}}c}
\includegraphics[height=1.3in]{sta20061026figure.3} &
\includegraphics[height=1.3in]{sta20061026figure.7}&
\includegraphics[height=1.3in]{sta20061026figure.6} \\
\multicolumn{1}{c}{\mbox{\bf (a)}} &
    \multicolumn{1}{c}{\mbox{\bf (b)}} &
    \multicolumn{1}{c}{\mbox{\bf (c)}}\\
    \includegraphics[height=1.3in]{sta20061026figure.4} &
\includegraphics[height=1.3in]{sta20061026figure.5}&
\includegraphics[height=1.3in]{sta20061026figure.8} \\
\multicolumn{1}{c}{\mbox{\bf (d)}} &
    \multicolumn{1}{c}{\mbox{\bf (e)}} &
    \multicolumn{1}{c}{\mbox{\bf (f)}}
\ea$
\ec
\caption{The Riemann problem for the LWR model: Stationary states in supply-demand diagrams}\label{fig:type1} \efg

For the original LWR model \refe{lwr}, the upstream link 1 and the downstream link 2 have the same fundamental diagram.
Therefore we have $C_1=C_2$.
In \citep{lebacque1996godunov}, there are four scenarios for solutions to the Riemann problem. Here we re-organize them into the following six cases of initial conditions:
\ben
\item [Case 1] Link 1 is SUC, and link 2 is UC. That is, $D_1< S_1=C_1$ and $D_2\leq S_2=C_2$. In this case, $D_1 < S_2$.  From Theorem \ref{linearboundary} we have that $U_1^-=U_1(0^-,t)=U_2^+=U_2(0^+,t)=U_1$ and $q_{1\to2}(0,t)=q(U_1)$. Therefore, on link 1, there is no wave; and on link 2, there is a forward shock wave when $D_1<D_2$, a forward rarefaction wave when $D_1>D_2$, and no wave when $D_1=D_2$.
\item [Case 2] Link 1 is OC, and link 2 is UC. That is, $S_1 \leq D_1=C_1$ and $D_2 \leq S_2=C_2$. In this case, $D_1=S_2=C_2$.  From Theorem \ref{linearboundary} we have that $U_1^-=U_1(0^-,t)=U_2^+=U_2(0^+,t)=(C_1,C_2)$  and $q_{1\to2}(0,t)=C_1$. Therefore, on link 1, there is a backward rarefaction wave when $S_1<C_1$ and no wave when $S_1=C_1$; and on link 2, there is a forward rarefaction wave when $D_2<S_2$ and no wave when $D_2=S_2$.
\item [Case 3] Link 1 is OC, and link 2 is SOC. That is, $S_1 \leq D_1=C_1$ and $S_2 < D_2=C_2$. In this case, $C_1=D_1 > S_2$.  From Theorem \ref{linearboundary} we have that $U_1^-=U_1(0^-,t)=U_2^+=U_2(0^+,t)=U_2$  and $q_{1\to2}(0,t)=q(U_2)$. Therefore, on link 1 there is a backward shock wave when $S_1>S_2$, a backward rarefaction wave when $S_1<S_2$, and no wave when $S_1=S_2$; and on link 2, there is no wave.
\item [Case 4] Link 1 is SUC, and link 2 is OC, and $q(U_1)<q(U_2)$. That is, $D_1< S_2\leq D_2=S_1=C_1$. From Theorem \ref{linearboundary} we have that $U_1^-=U_1(0^-,t)=U_2^+=U_2(0^+,t)=U_1$  and $q_{1\to2}(0,t)=q(U_1)$. Therefore, on link 1, there is no wave; and on link 2, there is a forward shock wave.
\item [Case 5] Link 1 is SUC, and link 2 is SOC, and $q(U_1)>q(U_2)$. That is, $S_2< D_1\leq S_1=D_2=C_1$. From Theorem \ref{linearboundary} we have that $U_1^-=U_1(0^-,t)=U_2^+=U_2(0^+,t)=U_2$   and $q_{1\to2}(0,t)=q(U_2)$. Therefore, on link 1, there is a backward shock wave; and on link 2, there is a no wave.
\item [Case 6] Link 1 is SUC, and link 2 is SOC, and $q(U_1)=q(U_2)$. That is, $D_1= S_2< D_2=S_1=C_1$. From Theorem \ref{linearboundary} we have that $U_1^-=U_2^+=(D_1,S_2)$, and $q_{1\to2}(0,t)=q(U_1)=q(U_2)$. Therefore, on link 1, there is no wave; and on link 2, there is no wave. In this case, there can exist interior states on link 1 or link 2: $U_1(0^-,t)=U_1^-$ and $\min \{ D_2(0^+,t), S_2(0^+,t)\} \geq D_1$ or $U_2(0^+,t)=U_2^+$ and $\min \{ D_1(0^-,t), S_1(0^-,t)\} \geq D_1$.
\een
Obviously the results of stationary states and kinematic waves above are consistent with those in \citep{lebacque1996godunov}. That is, the new solution framework yields the same wave solutions as traditional approaches.

The solutions of stationary states for the six cases are also shown in \reff{fig:type1}, where figures (a)-(f) are for cases 1-6 respectively. In these figures, both the upstream and downstream links share the same supply-demand diagram. From initial conditions $U_1$ and $U_2$ we can first draw the pair $(D_1,S_2)$, from which we can determine upstream and downstream stationary states accordingly.
Further we can summarize the solutions of stationary states in \reff{fig:lwr-sta} in the $(D_1, S_2)$ space.
This figure also demonstrates a graphical scheme for solving the stationary states as follows. First, from initial $U_1$ we draw a vertical line (thin pink line with an arrow),  from initial $U_2$ we draw a horizontal line (thin pink line with an arrow), and the intersection point is $(D_1, S_2)$. Then, if the intersection point is above the line $0A$, we draw a vertical line (thick blue line with arrow), and its intersection with $AC$ gives the stationary states; if the intersection point is below the line $0A$, we draw a horizontal line (thick blue line with arrow), and its intersection with $AC$ gives the stationary states; if the intersection point is on the line $0A$, we draw both a vertical line (thick blue line with arrow) and a horizontal line (thick blue line with arrow), and their intersections with $AC$ are the stationary states for the upstream and the downstream links respectively.
Note that this scheme also works when the upstream and downstream links do not have the same fundamental diagrams but the same supply-demand diagram, i.e., $C_1=C_2$.

\bfg\bc
\includegraphics{sta20061026figure.9}\ec
\caption{Solution of stationary states for the Riemann problem for the LWR model}\label{fig:lwr-sta}
\efg

\subsection{The inhomogeneous LWR model}
When $C_1\neq C_2$, then the road is inhomogeneous, and there is a discontinuity in the fundamental diagram at $x=0$.
In \reff{fig:inhlwr-sta}, we demonstrate a graphical scheme for solving stationary states in the $(D_1, S_2)$ space for the inhomogeneous LWR model.
We take \reff{fig:inhlwr-sta}(a) as an example, in which $C_1<C_2$.
 First, from initial $U_1$ we draw a vertical line (thin pink line with an arrow), from initial $U_2$ we draw a horizontal line (thin pink line with an arrow), and the intersection point is $(D_1, S_2)$. Then, if the intersection point is above the line $0A$, we draw a vertical line (thick blue line with arrow), and its intersections with $AC_1$  and $AC_2$ are the stationary states on links 1 and 2 respectively; if the intersection point $(D_1, S_2)$ is below the line $0A$, we draw a horizontal line (thick blue line with arrow), and its intersections with $AC_1$ and $AC_2$ are the stationary states on links 1 and 2 respectively; if the intersection point $(D_1, S_2)$ is on the line $0A$, we draw both a vertical line (thick blue line with arrow) and a horizontal line (thick blue line with arrow), and their intersections with $AC_1$ and $AC_2$ are the stationary states for the upstream and the downstream links respectively.
This scheme is the same as that for the homogeneous LWR model.

\bfg\bc $\ba{c@{\hspace{0.3in}}c}
\includegraphics[height=2in]{sta20061026figure.10} &
\includegraphics[height=2in]{sta20061026figure.11} \\
\multicolumn{1}{c}{\mbox{\bf (a)}} &
    \multicolumn{1}{c}{\mbox{\bf (b)}}
\ea$ \ec \caption{Solution of stationary states of the Riemann problem for the inhomogeneous LWR model}\label{fig:inhlwr-sta} \efg

In \citep{jin2003inhLWR}, the Riemann problem for the inhomogeneous LWR model was solved as a resonant nonlinear system, and ten types of wave solutions were obtained.
For example, wave solutions of Type 1 can be obtained in the new framework as follows. Both $U_1$ and $U_2$ are UC, $D_2<D_1\leq C_2=S_2$, and $C_1$ may be greater or smaller than $C_2$. From \reff{fig:inhlwr-sta} or Theorem \ref{linearboundary}, we can see that $U_1^-=(D_1, C_1)=U_1$, $U_2^+=(D_1, C_2)$, there is no wave on link 1, there is a forward rarefaction wave on link 2, and $q_{1\to 2}(0,t)=q(U_1)$.
It is easy to check that the wave solutions of other types are also consistent.

\section{Asymptotic traffic dynamics on an inhomogeneous ring road}
In this section we consider the inhomogeneous ring road with length $L$ shown in \reff{ringroad}, in which the traffic direction is shown by the arrow. The ring road is composed of two homogeneous links: link 1 with capacity $C_1$ for $x\in[0,L_1]$, link 2 with capacity $C_2$ for $x\in[L_1,L]$, the upstream boundary of link 1 is denoted as boundary 1, and the downstream boundary as boundary 2. Here we assume that link 1 is a bottleneck; i.e., $C_1<C_2$. For example, such a bottleneck can be caused by a smaller number of lanes. We assume the fundamental relationships for two links as $q=Q_1(\r)$ and $\r=R_1(\gamma)$ for $x\in[0,L_1]$, and $q=Q_2(\r)$ and $\r=R_2(\gamma)$ for $x\in[L_1,L]$.

\bfg\bc
\includegraphics[height=2in]{sd20080619figure.19}\ec
\caption{A ring road}\label{ringroad}
\efg

\subsection{Asymptotic stationary and interior states}
When the ring road reaches asymptotic stationary states, the flux at any location is the same, e.g., $q$. As we know, the asymptotic stationary state on a link can be uniformly UC, uniformly SOC, or a stationary shock wave (SS) connecting an upstream SUC state and a downstream SOC state \citep{bultelle1998shock}.  Then all possible combinations of stationary states are listed in \reft{ringsta} and explained in the following.
\bi
\item When link 1 is UC at $(q, C_1)$ with $q \leq C_1$, we have the following scenarios. (a) From Theorem \ref{linearboundary} it is possible that link 2 is UC at $(q,C_2)$, and the total number of vehicles on the ring road is $N_a=R_1(q/C_1) L_1+R_2(q/C_2) (L-L_1)$. (b) If link 2 is SS with upstream $(q,C_2)$ and downstream $(C_2,q)$, we have that $q=C_1$ and link 1 is critical at $(C_1,C_1)$. Assuming that link 2 is SUC for $x\in[L_1,L_2]$ and SOC for $x\in[L_2,L]$. In this case, the total number of vehicles on the ring road is $N_b=R_1(1) L_1+R_2(C_1/C_2) (L_2-L_1)+R_2(C_2/C_1) (L-L_2)$. (c) If link 2 is SOC at $(C_2,q)$,  we have from Theorem \ref{linearboundary} that $q=S_2=C_1$ at boundary 1. That is, link 1 is critical at $(C_1,C_1)$. In this case, the total number of vehicles on the ring road is $N_c=R_1(1) L_1+R_2(C_2/C_1) (L-L_1)$.
\item When link 1 is SOC at $(C_1,q)$ with $q<C_1$, we have the following scenarios. (d) It is possible that link 2 is SOC at $(C_2,q)$, and the total number of vehicles on the ring road is $N_d=R_1(C_1/q)L_1+R_2(C_2/q) (L-L_1)$.
    If link 2 is UC at $(q, C_2)$, we have from Theorem \ref{linearboundary} that $q=C_1$ at boundary 2. If link 2 is SS with upstream $(q,C_2)$, we have from Theorem \ref{linearboundary} that $q=C_1$ at boundary 2. Thus these two scenarios are impossible, since $q=C_1$ contradicts $q<C_1$.
\item When link 1 is SS with upstream $(q,C_1)$ and downstream $(C_1,q)$ with $q<C_1$, we have the following scenarios. If link 2 is UC at $(q,C_2)$, we have from Theorem \ref{linearboundary} that $q=\min\{C_1,C_2\}=C_1$ at boundary 2; if link 2 is SS with upstream $(q,C_2)$, we have from Theorem \ref{linearboundary} that $q=\min\{C_1,C_2\}=C_1$ at boundary 2; if link 2 is SOC at $(C_2,q)$, we have from Theorem \ref{linearboundary} that $q=\min\{C_2,C_1\}=C_1$ at boundary 1. All three of these scenarios are impossible, since $q=C_1$ contradicts $q < C_1$.
\ei

\btb\bc
\begin{tabular}{|c|c|c|c|}\hline
\backslashbox{Link 1}{Link 2}& UC $(q, C_2)$ & SS  $(q, C_2)\to (C_2, q)$ & SOC $(C_2, q)$\\\hline
UC $(q,C_1)$ & (a) & (b) & (c) \\\hline
SOC $(C_1,q)$ &x & x& (d)\\\hline
SS  $(q, C_1)\to (C_1, q)$& x& x& x\\\hline
\end{tabular}
\caption{All possible stationary states on the ring road in \reff{ringroad}} \label{ringsta}
\ec
\etb

For all four scenarios of asymptotic stationary states on the ring road, different scenarios have a different number of vehicles since
\bqs
0\leq N_a\leq R_1(1) L_1+R_2(C_1/C_2) (L-L_1) < N_b  <N_c=\\r_i(1) L_1+R_2(C_2/C_1) (L-L_1)<N_d\leq R_1(\infty) L_1+R_2(\infty) (L-L_1).
\eqs
Due to traffic conservation on the ring road, we can therefore determine the final stationary states by the initial number of vehicles $N$ on the road as follows: (a) When $N\leq R_1(1) L_1+R_2(C_1/C_2) (L-L_1)$, links 1 and 2 will be asymptotically stationary at UC with $(q,C_1)$ and $(q,C_2)$ respectively, where $q$ is the solution of $R_1(q/C_1) L_1+R_2(q/C_2) (L-L_1)=N$; (b) When $R_1(1) L_1+R_2(C_1/C_2) (L-L_1) < N<R_1(1) L_1+R_2(C_2/C_1) (L-L_1)$, link 1 will be asymptotically stationary at critical with $(C_1,C_1)$, and link 2 at SS with $(C_1,C_2)$ for $x\in[L_1, L_2]$ and $(C_2,C_1)$ for $x\in(L_2, L]$, where $L_2$ is the solution of $R_1(1) L_1+R_2(C_1/C_2) (L_2-L_1)+R_2(C_2/C_1) (L-L_2)=N$; (c) When $N=R_1(1) L_1+R_2(C_2/C_1) (L-L_1)$, link 1 will be asymptotically stationary at critical with $(C_1,C_1)$, and link 2 at SOC with $(C_2,C_1)$; (d) When $N>R_1(1) L_1+R_2(C_2/C_1) (L-L_1)$, links 1 and 2 will be asymptotically stationary at SOC with $(C_1,q)$ and $(C_2,q)$ respectively, where $q$ is the solution of $R_1(C_1/q)L_1+R_2(C_2/q) (L-L_1)=N$.

From Theorem \ref{linearboundary}, an interior state can occur at a boundary when its upstream demand equals the downstream supply, and its flux cannot be smaller than the demand or supply. In the following we consider possible asymptotic interior states on the ring road in \reff{ringroad}. First, at any location inside a uniform traffic stream on a homogeneous road, it is not possible to have interior states, since the upstream and downstream states are exactly the same at $(D, S)$ and $D=S$ if and only if the traffic is $D=S=C$, in which case the interior states have to be the same as the stationary states. Thus, interior states can only exist around the interface between two uniform traffic streams when the upstream demand equals the downstream supply, and we examine possible interior states in all four scenarios as follows. (a) The necessary condition for an interior state to exist at boundary 1 is $q=C_1$, i.e., when $N=R_1(1) L_1+R_2(C_1/C_2) (L-L_1)$. From Theorem \ref{linearboundary}, the interior state can only exist at $x=0^-$, but not $x=0^+$. The necessary condition for an interior state to exist at boundary 2 is $q=C_2$, which is not possible. (b) It is not possible for interior states to exist at either boundary 1 or 2, but it is possible for an interior state to exist around the SS interface at $x=L_2$. From Theorem \ref{linearboundary}, the interior state can only exist at $x=L_2^-$ or $x=L_2^+$. (c) It is not possible for an interior state to exist at boundary 1, but it is possible for an interior state to exist around boundary 2. That is, when $N=R_1(1) L_1+R_2(C_1/C_2) (L-L_1)$, from Theorem \ref{linearboundary}, the interior state can only exist at $x=L_1^+$, but not $x=L_1^-$. (d) It is not possible for  interior states to exist at either boundary 1 or 2.

In summary, there can exist three types of interior states: (a) When $N=N_1=R_1(1) L_1+R_2(C_1/C_2) (L-L_1)$, an interior state can only exist at $x=0^-=L^-$; (b) When $R_1(1) L_1+R_2(C_1/C_2) (L-L_1)<N<R_1(1) L_1+R_2(C_2/C_1) (L-L_1)$, an interior state can exist at $x=L_2^-$ or $x=L_2^+$; (c) When $N=N_3=R_1(1) L_1+R_2(C_2/C_1) (L-L_1)$, an interior state can only exist at $x=L_1^+$.

\subsection{Numerical examples}
In this subsection, we study asymptotic traffic dynamics on the inhomogeneous ring road in \reff{ringroad} with $L=600l$ = 16.8 km, $L_1=100l$=2.8 km, and the location-dependent speed-density relationships are based on \citep{kerner1994cluster,herrmann1998cluster} 
\bqs
V(\rho,a(x))=5.0461\left
[\left(1+\exp\{[\frac{\r}{a(x)\r_j}-0.25]/0.06\}\right)^{-1}-3.72\times
10^{-6}\right] l/\tau,
\eqs
where the relaxation time $\tau$ = 5 s;  the unit length $l$ = 0.028 km; the free flow speed $v_f= 27.8$ m/s; the jam density of a single lane $\r_j$ = 180 veh/km/lane. Here the number of lanes $a(x)=1$ for link 1 and $a(x)=2$ for link 2.
The corresponding fundamental diagram $q=Q(\r, a(x))$ is non-convex but unimodal in density $\r$.
In addition, $C_1=0.7091$ veh/s, and $C_2=2C_1$.
Thus we can compute $R_1(1)$=35.8944 veh/km, $R_2(C_1/C_2)=R_2(\frac 12)$=26.4162 veh/km, and $R_2(C_2/C_1)=R_2(2)$=118.3550 veh/km. Hence $N_1=R_1(1) L_1+R_2(C_1/C_2) (L-L_1)=470.3311$ veh, and $N_3=R_1(1) L_1+R_2(C_2/C_1) (L-L_1)=1757.4746$ veh.

Here we consider the following initial condition:
\bqn
\ba{ccll}
\rho(x,0)&=&a(x)(\r_0+3 \sin \frac {2\pi x}L), &x\in[0,L],
\\
v(x,0)&=&V(\r(x,0),a(x)),&x\in[0,L].
\ea\label{global2}
\eqn
Then, the total number of vehicles on the ring road is
\bqs
N&=&2\r_0L-\int_{0}^{100l} (\r_0+3 \sin \frac {2\pi x}L) dx=1100 l \r_0-\frac{450}{\pi} l
\eqs
When $\rho_0=15.4007$ veh/km, $N=N_1$ and we observe an interior state at $x=0^-=L^-$; when $\rho_0=57.1911$ veh/km, $N=N_3$ and we observe an interior state at $x=L_1^+$; when $\rho_0\in(15.4007,57.1911)$, we observe an interior state at $x=L_2^-$ or $x=L_2^+$, where $L_2$ is the solution of $R_1(1) L_1+R_2(\frac12) (L_2-L_1)+R_2(2) (L-L_2)=N$. For example, when $\rho_0=28$ veh/km, $N=858.3893$ veh, and
\bqs
L_2&=&\frac{N-(R_1(1) -R_2(\frac12)) L_1-R_2(2)L}{R_2(\frac12)-R_2(2)}=449.2561l.
\eqs

\bfg\bc\includegraphics[height=3in]{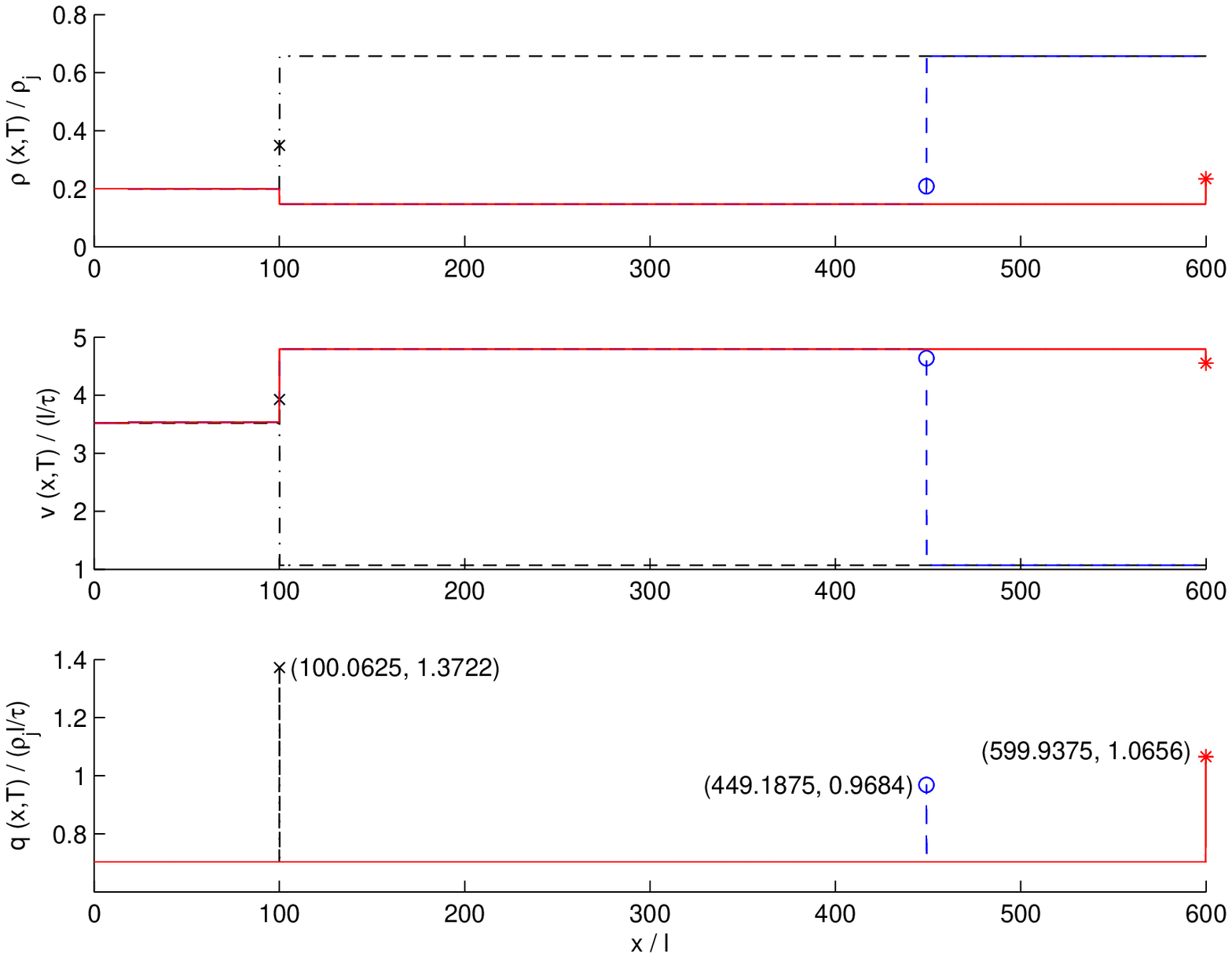}\caption{Solutions of $\r$, $v$, and $q$ at $T=24000$ s for initial conditions in \refe{global2}: solid lines with stars for $\r_0=15.4007$ veh/km, dashed lines with circles for $\rho_0=28$ veh/km, and dash-dotted lines for $\rho_0=57.1911$.} \label{stationary20080709plot}
\ec\efg

In the following, we simulate traffic dynamics on the ring road for three different initial $\r_0$: 15.4007 veh/km, 28 veh/km, and 57.1911 veh/km. Here we use the Godunov finite difference equation in \refe{lwrfd} and the supply-demand method in \refe{sdmethod} for computing boundary fluxes. The simulation time is  $T=4800 \tau$ = 24000 s. We partition the road $[0, L]$ into $N=4800$ cells and the time interval $[0, T]$ into $K=240000$ steps. Hence, the length of each cell is $\dx=3.5$ m and the length of each time step is $\dt=0.1$ s. The CFL condition number \citep{courant1928CFL} is $v_f \frac {\dt}{\dx}\leq 0.79<1$. The results for the three initial conditions are shown in \reff{stationary20080709plot}, where the bottom figure shows the locations and fluxes of all three interior states. From the figure, we can see that each of the three interior states only exists in one cell, and the locations of interior states are exactly as predicted above. Note that the top right figure of Figure 18 in \citep{jin2003inhLWR} also demonstrates the existence of an interior state, which is at the interface of a stationary shock.

\section{Conclusion}
In this paper we first reviewed the definitions of the supply and demand functions and the discrete supply-demand method for computing boundary fluxes. We then introduced the supply-demand diagram of a roadway and a new framework for solving the Riemann problem of the inhomogeneous LWR model in supply-demand space. In this framework, each link can have asymptotic interior and stationary states near the boundary, and the wave on each link is determined by the Riemann problem of the homogeneous LWR model with stationary and initial states for initial conditions.
We have derived conditions for admissible stationary and interior states and introduced an entropy condition based on the discrete supply-demand method for computing boundary fluxes.
We then proved that solutions to the Riemann problem exist and are unique and demonstrated that these solutions are consistent with those in literature for both the homogeneous and inhomogeneous LWR models. We also presented a graphical approach for finding the asymptotic stationary states with the help of supply-demand diagrams. Finally, we discussed the asymptotic stationary states on a ring road with arbitrary initial conditions and demonstrated with numerical examples that the existence and properties of the interior states are as predicted in this framework.

Unlike existing studies of the homogeneous or inhomogeneous LWR models, this study analyzes traffic dynamics in supply-demand space. In this framework, the discrete supply-demand method is applied as an entropy condition. In this sense, our study provides a new approach for constructing convergent solutions of finite difference equations arising in a Godunov method \refe{lwrfd} with a supply-demand method \refe{sdmethod} for both the homogeneous and inhomogeneous LWR models.
We have demonstrated that this new approach can successfully predict the existence and properties of interior states in numerical solutions.  However, note that interior states take only one cell in numerical solutions and vanish as we diminish the cell size. In this sense, the interior states are inconsequential to solutions of the Riemann problem.

Compared with existing studies, the new approach in the supply-demand framework is much simpler. In addition, since supply-demand methods have been proposed for computing fluxes through other junctions in general road networks \citep{daganzo1995ctm,lebacque1996godunov,jin2003merge,jin2003dissertation}, our framework could be extended to constructing solutions to the Riemann problem in these models.
In \citep{jin2009merge}, we successfully applied this framework to analyze the Riemann problem of merging traffic flow. In addition, one could also apply this new framework to analyze asymptotic traffic dynamics in a road network, such as the diverge-merge network studied in \citep{jin2008_network}.

\begin{acknowledgement}
The work of the first two authors was supported in part by the National Natural Science Foundation of China (No. 50708107), the Hi-Tech Research and Development Program of China (863 Project) (No. 2007AA11Z222), and the National Basic Research Program of China (973 Project) (No. 2006CB705506). The work of the third author was supported in part by the US Department of Energy (DOE) Mathematical, Information, and Computing Sciences Division under contract number DE-FG02-03ER25579.
\end{acknowledgement}

\end {document}